\documentclass[11pt]{article}
\usepackage[french, english]{babel}
\usepackage[utf8]{inputenc}
\usepackage{amsmath,amsfonts,amssymb,amsthm,mathrsfs}
\usepackage{minitoc}
\usepackage{bbm}
\renewcommand{\leq}{\leqslant}
\renewcommand{\geq}{\geqslant}
\usepackage[mono=false]{libertine}
\useosf
\usepackage{mathpazo}

\usepackage{wasysym}

\usepackage[dvipsnames]{xcolor}
\usepackage[a4paper,vmargin={3.5cm,3.5cm},hmargin={2cm,2cm}]{geometry}
\usepackage[font=sf, labelfont={sf,bf}, margin=1cm]{caption}
\usepackage{graphicx,graphics}
\usepackage{epsfig}
\usepackage{latexsym}
\usepackage{stmaryrd}
\usepackage{ae,aecompl}

\usepackage[english]{babel}
 \usepackage[colorlinks=true]{hyperref}
\usepackage{pstricks}
\usepackage{enumerate}
\usepackage{tikz}
\usepackage{fontawesome}
\usetikzlibrary{arrows,automata}

\newcommand{\E}{\mathbb{E}}
\DeclareFixedFont{\beaupetit}{T1}{ftp}{b}{n}{2cm}

\newtheorem{theorem}{Theorem}[]
\newtheorem{definition}{Definition}[]
\newtheorem{proposition}[]{Proposition}

\newtheorem{lemma}[]{Lemma}
\newtheorem{corollary}[]{Corollary}

\theoremstyle{definition}

\usepackage{todonotes}

\title{\textsc{Parking on trees with a (random) given degree sequence and the Frozen configuration model}}
\author{
Alice \textsc{Contat}\thanks{Universit\'e Sorbonne Paris-Nord.\hfill  \href{mailto:alice.contat@math.cnrs.fr}{\texttt{alice.contat@math.cnrs.fr}}}
}

\date{}
\begin{document}
\maketitle 

\begin{abstract} Consider a rooted tree $ \mathfrak{t}$ on the top of which we let cars arrive on its vertices. Each car tries to park on its arriving vertex but if it is already occupied,  it drives towards the root of the tree and parks as soon as possible. In this article, we establish a natural coupling between the parking process on trees with prescribed degrees and an oriented configuration model. As a consequence, we recover the location of the phase transition for parking on critical Bienaym\'e--Galton--Watson trees already proven in \cite{contat2020sharpness,curien2022phase}.
\end{abstract}

\section{Introduction}

Parking problems on random trees received a lot of attention in the last decade. Initially introduced by Konheim and Weiss in the case of a line \cite{konheim1966occupancy}, their study on trees has been initiated by Lackner and Panhozler in 2016 \cite{LaP16}. Since then, an intriguing phase transition has been brought to light especially on critical Bienaym\'e--Galton--Watson trees \cite{chen2021parking,contat2020sharpness,contat2021parking,curien2022phase,GP19, JO18} and on the infinite binary tree \cite{aldous2022parking}. The goal of this work is to shed a new light on this phase transition already observed in \cite{contat2020sharpness} by coupling the parking model with a random graph model.   The apparition of a macroscopic flux of outgoing for the parking process corresponds to the apparition of a giant component for our graph model. A similar coupling already appears in  \cite{contat2021parking} for a specific type of trees and specific car arrivals law. We extend here this technique for a very  general context and we plan to use this tool to study the critical behavior for a large class of parking models. 

\paragraph{Parking on trees.} Let us be more precise on the parking process and consider a (fixed) finite rooted tree $\mathfrak{t}$. This tree represents the parking lot where each vertex can accommodate at most one car and we imagine that all edges are oriented towards the root of $ \mathfrak{t}$. On the top of $ \mathfrak{t}$, we put a (fixed) car arrival decoration $ (a_x : x \in \mathfrak{t})$. 
The parking rule is the following. Each car tries to park on its arriving vertex and if it is already occupied, it drives towards the root, following the edges and take the first available spot. If no free spot is found during its descent to the root, then the car exits the tree and contribute to the \emph{flux} of outgoing cars. We write $ \varphi ( \mathfrak{t)}$ the \emph{flux} of outgoing cars on $ \mathfrak{t}$. \\
An important property of this model is its Abelian property: the final configuration of occupied spots and the number of outgoing cars do not depend on the order in which the cars park.
Note also that the parking process only depends on the geometry of the underlying tree but not on the potential labels of the vertices. In what follows, we will however consider graphs with labeled vertices to avoid symmetry issues. \\
Our purpose is to study the parking process not on fixed but on random trees with a fixed degree sequence. The goal of this work is to link the parking process on those random trees to a variation around the configuration model. Our coupling is the basis to derive numerous properties of the parking process such that its phase transition and its location, and explore its different universality classes. 

\paragraph{Random  trees with given degrees.} 
We start by describing the model of trees which we consider. We fix a sequence $ \mathbf{i}=(i_1, \dots, i_n)$  such that $ \sum_{k \geq 1} i_k = n-1$, so that $i_k$ is the number of children of the vertex with label $k$. 
Then, the number of \emph{(rooted) plane trees} with labeled vertices such that the vertex label $k$ has $i_k$ children for all $k \geq 1$ is 
$$ (n-1)!,$$
see for example \cite{harary1964number}, and we write $ T ( \mathbf{i})$ such a tree chosen uniformly at random.
We will couple the parking process on $ T ( \mathbf{i})$ with a variation of the configuration model which we now present.

\paragraph{Frozen configuration model.} From now on, we consider oriented graphs in a broad sense since we allow vertices to have unmatched (oriented!) half-edges or \emph{legs}. We fix $n \geq 1$ and fix $ \mathbf{i} = (i_1, \dots, i_n)$ an in-degree sequence and $ \mathbf{o} = (o_1, \dots, o_n)$ an out-degree sequence, which means that the vertex with label $k$ has in-degree $i_k$ and out-degree $o_k$. The oriented configuration model $ \mathrm{CM} ( \mathbf{i}, \mathbf{o})$ is the graph obtained by starting with $n$ labeled vertices with this degree sequences and sampled by pairing as many in- and out-legs as possible. An important property is that we can couple the legs step-by-step. Indeed, suppose that $ S_ \mathrm{in} \leq S_ \mathrm{out}$ and start with vertices with the appropriate number of legs. At each step, we can choose an (unmatched) in-leg, where the choice only depends on the past, and match it with an unmatched out-leg sampled uniformly at random and independently from the past until there is no more unmatched in-leg. When $ S_ \mathrm{in} \geq S_ \mathrm{out}$, the same construction is valid when exchanging the role of in- and out-legs.

Note that unlike most of the existing literature, we allow  the total out-degree $ S_{ \mathrm{out}} := \sum_{k \geq 1 } o_k $ to be different from the total in-degree $ S_{ \mathrm{in}} := \sum_{k \geq 1 } i_k$, so that they may be unmatched in-legs or out-legs (but not both) in  $ \mathrm{CM} ( \mathbf{i}, \mathbf{o})$. 

\begin{figure}[!h]
 \begin{center}
 \includegraphics[width=17cm]{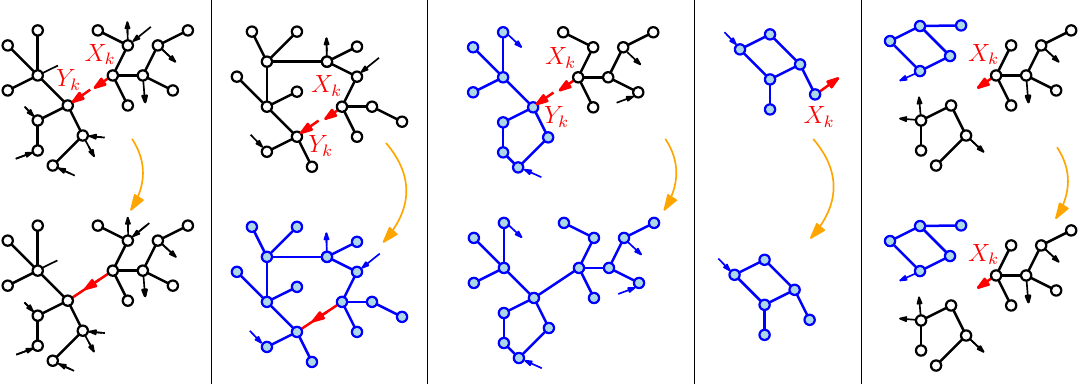}
 \caption{\label{fig:transitionsFCM} Illustration of the possible transition for the frozen configuration model. In the first three cases, the out-leg $X_k$ goes out from a with vertex, thus we sample an in-leg $Y_k$ uniformly at random and we match the two oriented legs to create an edge. On the left, the in-leg $ Y_k$ was incident to a white vertex and the addition of the edge does not create a cycle, thus the colors of the vertices remain unchanged. The two legs which we eventually match are displayed in red. In the second case, the addition of the edge creates a cycle and the whole component of the edge becomes frozen (blue). In the third case, the in-leg $Y_k$ was incident to a blue vertex and all vertices of the components of $X_k$ become blue. In the fourth case, the out-leg $X_k$ goes out from a blue vertex, so we delete it without sampling an in-leg $Y_k$. Finally, if there are no in-leg remaining, we delete the out-leg $X_k$.}
 \end{center}
 \end{figure}

We now introduce its \emph{frozen} version. We start from the same degree sequences and construct the graph step-by-step and slow down the growth of components which are not trees. The graph has oriented edges but we consider the \emph{weekly connected} components i.e.\ the connected components of the unoriented version of the graph (two vertices are in the same component if there exists an unoriented path of edges between them). 
 The vertices will be of two types, either white for ``standard" vertices or blue for ``frozen" vertices. We start from $n$ white vertices and try to match the out-legs one by one with the following rule, See Figure \ref{fig:transitionsFCM}: at each step, we choose an unmatched out-leg uniformly at random and independently of the past and
\begin{itemize}
\item if this out-leg goes out from a blue vertex, then we erase it without matching it.
\item Otherwise, if this out-leg goes out from a white vertex, we match it with an in-leg chosen uniformly at random (if there is one) and add this edge to the graph. If there is no more unmatched in-leg, we erase the out-leg.  If the vertex incident to the in-leg was blue or if the added edge creates a cycle or if there was no more unmatched in-leg, then we color all the vertices of the (weak) connected component in blue. 
\end{itemize}
Then the frozen configuration model is the graph $ \mathrm{FCM} ( \mathbf{i}, \mathbf{o})$ obtained when all unmatched out-legs are either matched or deleted. See Section \ref{sec:defFCM} for a more precise definition. Note that they may be unmatched in-legs in  $ \mathrm{FCM} ( \mathbf{i}, \mathbf{o})$ but no unmatched out-leg. This construction is the natural generalization of the \emph{frozen Erd\H{o}s--Rényi} random graph model introduced by Contat and Curien in \cite{contat2021parking}.

\paragraph{Coupling.} The main input of this paper is to establish a natural coupling between the frozen configuration model and the parking process on a tree with prescribed degree sequence and prescribed car arrivals. The main idea is to consider, on the parking side, the underlying tree as unknown and to park the cars one-by-one and construct the (random) tree step-by-step when we need to accommodate a car. An important observable of our model is the subforest spanned by the edges emanating directly from a vertex containing a parked car (recall that in case of trees, the edges are oriented towards the root). The connected components of this subforest are called the \emph{nearly parked components}.

We state here our main theorem in an informal way. See Section \ref{sec:couplingtree} for a precise statement. 

\begin{theorem}[informal] \label{thm:coupling}  Let us fix $n \geq 1$ and a in-degree sequence $ \mathbf{i} = (i_1, \dots, i_n)$ such that \linebreak[4] {$ S_{ \mathrm{in}}:=\sum_{k \geq 1} i_k = n-1$}, and a car decoration $ \mathbf{a} = (a_1, \dots, a_n)$.
We can couple the parking process on the random tree $T( \mathbf{i})$  with car arrivals prescribed by $ \mathbf{a}$ and the frozen configuration model $\mathrm{FCM}( \mathbf{i}, \mathbf{a})$ such that the nearly parked components of the parking process are the weakly connected components of the graph $\mathrm{FCM}( \mathbf{i}, \mathbf{a})$ except for the component of the root which contained all blue vertices. Moreover, the out-legs deleted while constructing $\mathrm{FCM}( \mathbf{i}, \mathbf{a})$ correspond to outgoing cars contributing to the flux.
\end{theorem}

This theorem is a natural generalization of \cite[Proposition 8]{contat2021parking}. As in \cite{contat2021parking}, the coupling is easier to understand in the case of parking on mappings and we start by presenting this coupling in Section  \ref{sec:couplingmapping}. The proof technique is essentially the same and consists in discovering or \emph{peeling} the underlying tree step-by-step instead of revealing the whole tree first.

\paragraph{Random degree sequences.} As an application of our theorem, we can have an additional layer of randomness and imagine that the degree sequences $ \mathbf{i}$ and $ \mathbf{a}$ are also random. In particular, we can encompass the model of Bienaym\'e--Galton--Watson conditioned to have size $n$ and recover the classification given by Curien and Hénard \cite[Theorem 1]{curien2022phase} and by Contat \cite[Theorem 1]{contat2020sharpness}. We work here in a slightly more general context. More precisely,  we consider a sequence of random degree sequences $ \mathbf{I}^{(n)}=  ({I}^{(n)}_1, \cdots, {I}^{(n)}_n)$ and $ \mathbf{A}^{(n)}=  ({A}^{(n)}_1, \cdots, {A}^{(n)}_n)$ such that for all $n$, the total in-degree $ \sum_{k=1}^{n } {I}_k^{(n)}$ is equal to $n-1$ and 
\begin{equation}\label{eq:hypo-convergence}\tag{$H_{ \mathrm{conv}}$}  \frac{1}{n} \sum_{k=1}^{n} \delta_{({I}^{(n)}_k, {A}^{(n)}_k)} \xrightarrow[n\to\infty]{} \lambda := \sum_{k \geq 0} \nu_k \sum_{j \geq 0} \mu_{(k),j} \delta_{(k,j)}, \end{equation}
where the measure $ \nu:=\sum_{k \geq 0} \nu_k \delta_k$ is a probability measure that we see as an offspring distribution and for all $k \geq 0$, the measure $ \mu_{(k)} =\sum_{j \geq 0} \mu_{(k),j} \delta_j$ is a probability measure that represents the typical distribution of the car arrivals on a vertex of out-degree $k$. We assume that $ \nu$ has mean $1$ and finite non zeri variance $ \Sigma^2 \in (0, \infty)$, and for all $k \geq 0$, we assume that $ \mu_{(k)} := \sum_{j \geq 0} \mu_{(k),j} \delta_j$ has mean $m_{(k)} < \infty$ and finite variance $ \sigma_{(k)}^2$.  All theses assumptions are refered as \eqref{eq:hypo-convergence} in the rest of the work. 

The parking model on random trees with a random degree sequence works as follows: we first sample the degree sequences $(\mathbf{I}^{(n)}, \mathbf{A}^{(n)})$ and we consider the random tree $  \mathrm{T}(\mathbf{I}^{(n)})$ such that conditionally on $\mathbf{I}^{(n)}$, it is a random tree with prescribed number of children given by $\mathbf{I}^{(n)}$. On the top on this tree, we have car arrivals given by $ \mathbf{A}^{(n)} = ( A_k^{(n)} : 1 \leq k \leq n)$ and park the cars with the usual parking rule. Recall that  $ \varphi(\mathrm{T} ( \mathbf{I}^{(n)}))$ denotes the flux of outgoing cars.

To state our second theorem, we need to introduce the  size-biased distribution $\overline{ \nu} = \sum_{k \geq 0} k \nu_{ k} \delta_k$ and the three quantities 
\[ \E_{\overline{\nu}} [m] := \sum_{k = 0}^{\infty} k \nu_k \, m_{(k)},  \quad \E_{\nu} [m] := \sum_{k = 0}^{\infty} \nu_k \, m_{(k)} \quad\text{ and } \quad \E_{\nu} [\sigma^2 + m^2 -m] := \sum_{k = 0}^{\infty}\nu_k  \left( \sigma_{(k)}^2 + m_{(k)}^2 -m_{(k)}\right).\] 

\begin{theorem} \label{thm:coupling:curienhenard}
We assume the hypothesis \eqref{eq:hypo-convergence} with $\E_{\overline{\nu}} [m] \leq 1$ and $\E_{{\nu}} [m] \leq 1$. We also assume that there exists a constant  $K$ such that $m_{(k)}< K$ and $ \sigma_{(k)}^2< K$ for all  $ k \geq 0$. The parking process undergoes a phase transition which depends on the sign of the quantity 
\begin{equation}\label{tTheta}
\Theta := (1- \E_{\overline{\nu}} [m])^2 -\Sigma^2 \E_{\nu} [\sigma^2 + m^2 -m].
\end{equation}
More precisely, we have 
\begin{equation*} 
\frac{ \varphi ( \mathrm{T} ( \mathbf{I}^{(n)}))}{n} \xrightarrow[n\to\infty]{( \mathbb{P})} C_{ \lambda}
\end{equation*}
 where $C_ { \lambda} = 0$ if and only if $ \Theta \geq 0$.
\end{theorem}

Our proof is based on the coupling with the random graph model given in Theorem \ref{thm:coupling}. An important observation is that the apparition of a macroscopic flux of outgoing cars ($ C_ \lambda >0$) coincides with the apparition of a giant component in our random graph model. Here we analyse this phase transition for the giant component by showing that the local limit of the graph is a (multitype) Bienaym\'e--Galton--Watson tree, which is almost surely finite when there is no giant component and that this coincides with the measures $ \lambda$ for which $ \Theta \geq 0$. 
This parameter $ \Theta$ is in fact the natural parameter which also appears in the physics literature by Kryven \cite{kryven2016emergence} for the apparition of the giant weakly connected component of the oriented configuration model. 

Note that the model considered by Contat in \cite{contat2020sharpness} and by Curien and Hénard in \cite{curien2022phase} falls into this framework. Indeed, if we consider a Bienaym\'e--Galton--Watson tree with offspring distribution $ \nu$ conditioned to have $n$ vertices, on which we add cars independently on each vertex and such that the law of the number of cars arriving on a vertex with $k$ children is $ \mu_{(k)}$, then its in-degree sequence and car arrival decoration satisfies the hypothesis \eqref{eq:hypo-convergence}, see for example \cite[Theorem 7.11]{Jan12b}.
However their proof is different. They use the magic of skip-free random walk to write a differential equation satisfied by the mean flux of outgoing cars (on unconditioned Bienaym\'e--Galton--Watson trees).

\paragraph{ Acknowledgments. } The author thanks Nicolas Curien for motivating discussions during the elaboration of this work.
This work has been supported by ANR RanTanPlan.

\section{Coupling between parking and oriented configuration model}

The goal of this section is to explain our main coupling between the parking process and the frozen configuration model, so that we can  formally state and prove Theorem \ref{thm:coupling}.
The coupling is easier to understand when the underlying graph is a random mapping (and not a tree) with a prescribed degree sequence. For the sake of clarity, we start by presenting this version. 
\subsection{Parking on mapping with given degree}\label{sec:parkingmapping}

\begin{figure}[!h]
 \begin{center}
 \includegraphics[width=14cm]{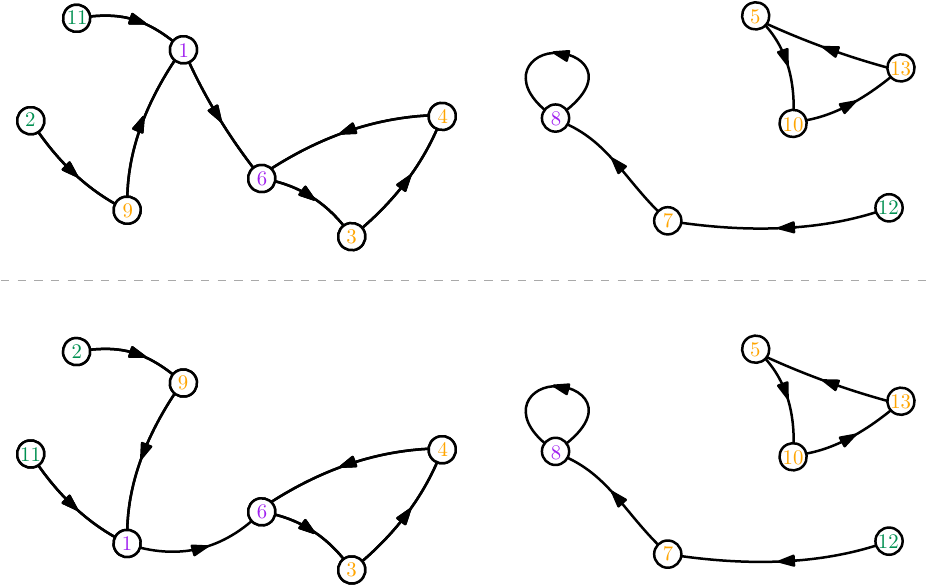}
 \caption{Two plane mappings with distinguishable vertices. The in-degree sequence is $ \mathbf{i} = (2,0,1,1,1,2,1,2,1,1,0,0,1)$ and the labels of the vertices of in-degree $0$, $1$ and $2$ are colored  green, orange and purple respectively. Note that these two mappings are different because the planar orderings around the vertex with (purple) label $1$ are not the same. \label{fig:mapping} }
 \end{center}
 \end{figure}

Recall that a \emph{mapping} is an oriented graph with labeled vertices, where there is exactly one edge going out from each vertex. Here, we will consider \emph{plane} mapping, which means that there is a cyclic orientation of the edges around each vertex, see Figure \ref{fig:mapping}. 
 
As in the case of tree, we can define a mapping with a prescribed degree sequence. We fix $n \geq 1$ and let $ \mathbf{i}=(i_1, \dots, i_n)$ be such that $ \sum_{k \geq 1} i_k = n$. For all $k \geq 1$, we imagine that the vertex with label $k$ has in-degree $i_k$ and out-degree $1$. All in-legs are distinguishable since there is a planar ordering of the legs of each vertex.  We consider the random (plane) mapping $ \mathrm{RM} ( \mathbf{i})$ with given degree $ \mathbf{i}$ obtained by pairing each out-leg with an in-leg uniformly at random. Note that the pairing is necessarily plane since there is at most one cycle per (weakly) connected component. It is easy to see that there are $ n!$ possible pairing of the legs. To see this, we can deterministically label the in- and out-legs, each from $1$ to $n$. For example, we can imagine that the out-leg of each vertex has the same label as the vertex and we label the in-legs vertex by vertex and clockwise for each vertex. Then a pairing of in- and out-edges is just a bijection between the labels of the in-legs and the labels of the out-legs. 
An important observation is that the pairing of the edges can be made in a progressive manner. Indeed we can choose one-by-one, deterministically or randomly, independently of the remaining pairing,  an unmatched out-leg, and for this out-leg choose an unmatched in-leg uniformly at random (and independently from the past) and match them together. We will see an illustration of this observation in Section \ref{sec:couplingmapping}.

Note that the labels of the vertices has an impact on the geometry of the mapping. More precisely, if for all $k \geq 0$, we write $d_k := \# \{j : i_j = k\}$ the number of vertices of in-degree $k$ in the sequence $ \mathbf{i}$. Then, the law of a mapping chosen uniformly at random among the (plane) mappings containing $d_k$ (unlabeled) vertices with $k$ children for all $k \geq 0$ is different from the law of the mapping obtained from $ \mathrm{RM} ( \mathbf{i})$ by forgetting the labels of its vertices. This is not the case for trees with prescribed number of children.

\paragraph{Parking on mappings.} As on trees, we want to apply the parking process on mappings. We fix a finite mapping $\mathfrak{m}$, and a car decoration $(a_x : x \in \mathfrak{m})$. The cars park in the same way as in trees: they follow the edges and take the first available spot. The only difference is that there is no root where the cars exit, but if a car is caught in an endless loop, then it exits the mapping without parking and contributes to the flux of outgoing cars, which we denote by $ \varphi ( \mathfrak{m})$. Another important notion is the submapping (subgraph of the mapping) of $ \mathfrak{m}$ spanned by all edges emanating from the occupied spots.  Its (weakly) connected components are called the \emph{near parked components} of $ \mathfrak{m}$.

We can also apply this parking rule to the random mapping $ \mathrm{RM} ( \mathbf{i})$ with a (fixed) car decoration $ \mathbf{a} = (a_k : 1 \leq k \leq n)$, and we write $ \varphi (\mathrm{RM} ( \mathbf{i}))$ the flux of outgoing cars.

\subsection{Frozen direct configuration model.}\label{sec:defFCM} 
We now define formally the frozen configuration graph model with which we will couple the parking process on random mappings. 
Let us start by recalling the usual oriented configuration model. We fix  $ \mathbf{i}=(i_1, \ldots, i_n)$ the in-degree sequence and  $ \mathbf{o}=(o_1, \dots, o_n)$ the out-degree sequence. To avoid issue about order between in- and out-edges and simplify the setting, we imagine that all in-legs are grouped together and  the out-legs are also grouped on another side (there is no alternance between out-legs and in-half edges) of each vertex. All legs are the distinguishable and there is a cyclic order of the legs around each vertex.  More precisely, we imagine each vertex has two sides (red and white) and that we put the out-legs on the white part and the  in-legs on the red part. With this convention, all legs are \emph{distinguishable} and there is a cyclic orientation of the legs around each vertex, as in the case of mapping with prescribed degree sequence. Unlike most of the existing literature, we allow  the total out-degree $ S_{ \mathrm{out}} := \sum_{k \geq 1 } o_k $ to be different from the total in-degree $ S_{ \mathrm{in}} := \sum_{k \geq 1 } i_k$.

\begin{figure}[!h]
 \begin{center}
 \includegraphics[width=8cm]{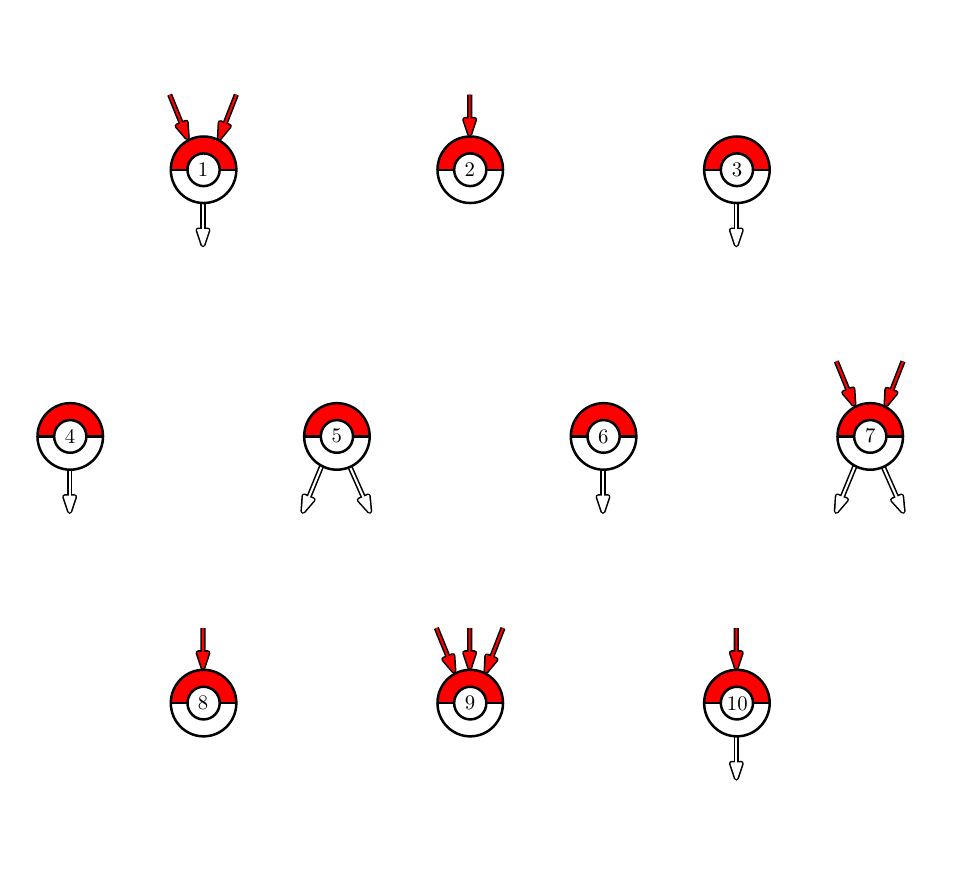} \hspace{0.5cm} \includegraphics[width=8cm]{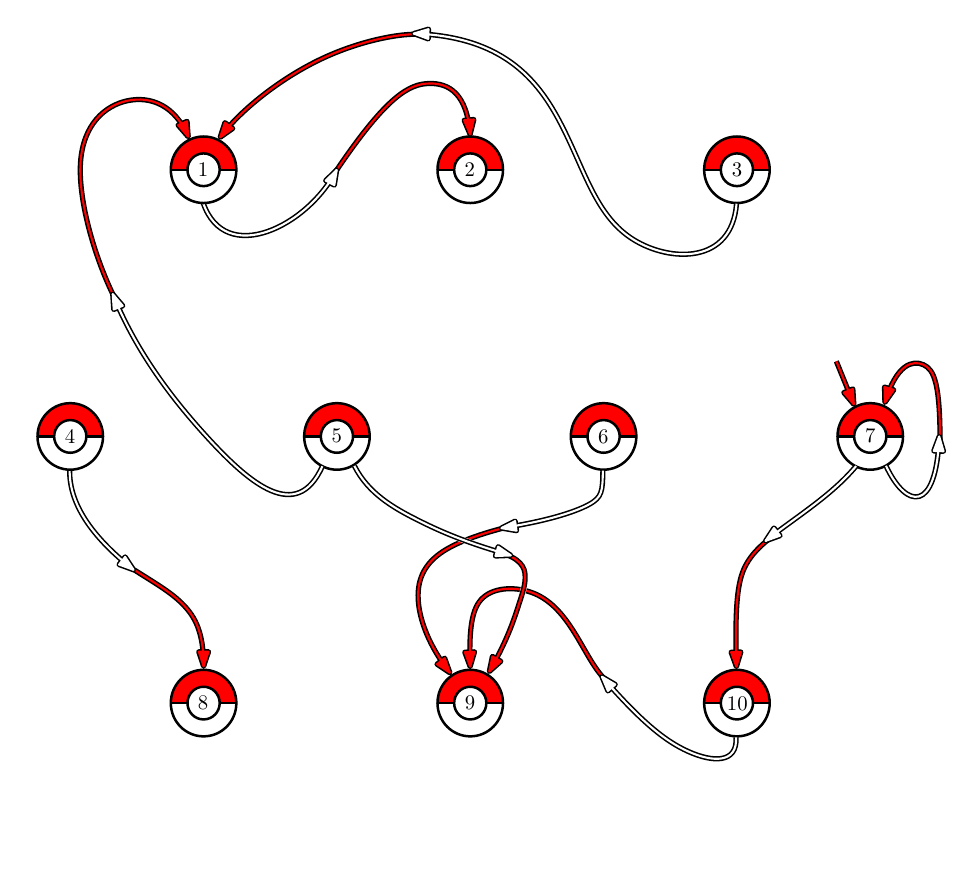}
 \caption{\label{fig:pokecm}Illustration of the general oriented configuration model. On the left, ten vertices with prescribed unmatched legs and on the right, a possible (maximal) pairing of the legs with one unmatched in-leg remaining.}
 \end{center}
 \end{figure}

We consider the random oriented multi-graph $ \mathrm{CM} ( \mathbf{i}, \mathbf{o})$ obtained by starting with $n$ labeled vertices where the vertex labeled $k$ has $i_k$ in-legs and  $o_k$ out-legs sampled by pairing the $S_{ \mathrm{out}}$ out-legs with (distinct) in-legs  if $S_{ \mathrm{out}} \leq S_{ \mathrm{in}}$ (resp.\ the $S_{ \mathrm{in}}$ in-legs with (distinct) out-legs  if $S_{ \mathrm{out}} > S_{ \mathrm{in}}$). Note that the number of edges in  $ \mathrm{CM} ( \mathbf{i}, \mathbf{o})$  is $ \min ( S_{ \mathrm{out}}, S_{ \mathrm{in}})$, and when $ S_{ \mathrm{out}} \neq S_{ \mathrm{in}}$, there are unmatched in- or out-legs remaining in $ \mathrm{CM} ( \mathbf{i}, \mathbf{o})$ (but not both), see Figure \ref{fig:pokecm}.
Since the legs are distinguishable, a nice property is that we can construct the pairing step by step. A way to see this step-by-step matching, when $ S_ \mathrm{out} < S_ \mathrm{in}$, is to start with vertices with the appropriate number of legs, and at each step, we  choose an (unmatched) in-leg, where the choice only depends on the past, and match it with an unmatched out-leg sampled uniformly at random and independently from the past until there is no more unmatched in-leg. Of course when $ S_ \mathrm{in} \geq S_ \mathrm{out}$, the same construction holds true when exchanging the role of in- and out-legs.

\paragraph{Frozen version.} Let us now introduce its frozen version. Again, we fix  $ \mathbf{i}=(i_1, \dots, i_n)$ the ``in-degree" sequence and the ``out-degree" sequence $ \mathbf{o}=(o_1, \dots, o_n)$ and for $1 \leq k \leq n$, we imagine that the vertex with label $k$ has two sides, a red side with $i_k$ out-legs and a white side with $o_k$ in-legs.  Again,  the total in-degree may be different from the total out-degree $ S_{ \mathrm{out}}  \neq S_{ \mathrm{in}}$, and this time, we will see that they may be unmatched in-legs remaining at the end of the step-by-step construction but no out-leg. We call the \emph{frozen configuration model} $ \mathrm{FCM} ( \mathbf{i}, \mathbf{o})$ the random graph obtained in a very similar way but slowing down a certain type of components. More precisely, the vertices of $ \mathrm{FCM} ( \mathbf{i}, \mathbf{o})$ are of two types, either white for ``normal" vertices or blue for the ``frozen" vertices. Note that these two colors live together with two colors for the two sides of each vertex. We start with $n$ white labeled vertices $\{1, \dots, n \}$ with the in-legs and out-legs described above. We then construct the graph sequentially: at each step $1 \leq k \leq S_{ \mathrm{out}} $, we sample 
\begin{center} \fbox{an unmatched out-leg $X_k$ uniformly at random and  independently from the past.} \end{center}
Sometimes, we delete this leg and eventually change the colors of the vertices with the following rules, See Figure \ref{fig:transitionsFCM}:
\begin{itemize}
\item If this out-leg goes out from a blue vertex, then we delete it without matching it.
\item Otherwise, if this out-leg goes out from a white vertex, we sample an unmatched in-leg $Y_k$ uniformly at random and  independently from the past (if there is one) and match the out-leg with the in-leg. If the in-leg was incident to a  blue vertex or if this  new edge creates a cycle or if there is no more unmatched in-half edge, then we color all the vertices of the (weak) connected component in blue. 
\end{itemize}

The graph obtained after step $S_{ \mathrm{out}}$ is $ \mathrm{FCM} ( \mathbf{i}, \mathbf{o})$. 
Note that they may be unmatched in-legs which we keep in the end.

\subsection{Coupling parking and Frozen configuration model.}\label{sec:couplingmapping}
Recall the notion of near components for the parking process, which are the connected components of the subgraph of the mapping spanned by the edges emanating from the occupied spots. Our first result is that we can couple the parking process on a mapping  with given degree and given car arrivals with the frozen configuration model. 
The following proposition is the analogue of Theorem \ref{thm:coupling} in the case of mappings. It is also similar to  \cite[Proposition 5]{contat2021parking} which treated the case of uniform (not planar) mapping without prescribed degrees and uniform car arrivals.

\begin{proposition} \label{prop:couplingmapping}Let us fix $ n \geq 1$, an in-degree sequence $ \mathbf{i}$ such that $ S_ \mathrm{in} := \sum_{k \geq 1} i_k =n$ and an out-degree sequence $ \mathbf{a}$. 
We can couple the parking process on the random mapping $ \mathrm{RM}( \mathbf{i})$  with car arrivals given by $ \mathbf{a}$ and the frozen configuration model $ \mathrm{FCM} ( \mathbf{i}, \mathbf{a})$ such that 
\begin{itemize}
\item The nearly parked components of  $ \mathrm{RM}( \mathbf{i})$ are the weakly connected component of $ \mathrm{FCM} (\mathbf{i}, \mathbf{a})$ (in terms of subset of vertices). 
\item The blue components in $ \mathrm{FCM} (\mathbf{i}, \mathbf{a})$ correspond to components which can acommodate no more cars in $ \mathrm{RM}( \mathbf{i})$.
\item the out-legs deleted while constructing $ \mathrm{FCM} ( \mathbf{i}, \mathbf{a})$ correspond to the cars  which did not manage to park in $ \mathrm{RM}(\mathbf{i})$.
\item  the mapping $ \mathrm{RM}( \mathbf{i})$ is obtained by pairing the remaining unmatched in-legs in a compatible way uniformly at random.
\end{itemize}

\end{proposition}

Note that in our construction of the configuration model, we chose to select the out-legs step-by-step uniformly at random. This means that we park the corresponding cars in this proposition uniformly at random. However, the parking process has an Abelian property: the final configuration does not depend on the order in which we park the cars. As a consequence, if we choose another algorithm, which only depends on the past,  to select the out-edge at each step (for example vertex-by-vertex with increasing labels, or the vertices with largest degree first...),  then the law of the sizes of weakly connected component are the same as that of the frozen configuration model. 

\begin{proof} The main idea is to construct the mapping step-by-step while we construct step-by-step the frozen configuration model using the same legs at each step.

\begin{center}
	\fbox{\begin{minipage}{15cm}
						
\noindent \textbf{From oriented edges to parking on mapping.} The coupling is made of two main parts. First, we start from a frozen configuration model $ \mathrm{FMC}( \mathbf{i}, \mathbf{a})$. Each out-leg is interpreted as a car arrival on a vertex and we will park them using the order given by the sequence $(X_k : 1 \leq k \leq S_{ \mathrm{out}})$ in the construction of the frozen configuration model (recall that by the Abelian property of the parking process, we can choose the order in which we want to park the cars). Since we use the same in-half-degree sequence $ \mathbf{i}$ for the mapping $ \mathrm{RM}(\mathbf{i})$, there is a natural correspondance between the in-legs of the mapping and that of the frozen configuration model.  If there is no ambiguity, we  use the same name for the two corresponding in-legs. In the mapping, we use the white sides of the vertices as the location of the edges going out from each vertex.

We use the out-legs (of the configuration model) and the in-legs as well as to construct iteratively an increasing sequence of ``submappings" $(\mathrm{RM}( \mathbf{i},\mathbf{a},m) : 0 \leq m \leq S_{ \mathrm{out}})$ where $ \mathrm{RM}( \mathbf{i}, \mathbf{a},0)$ is the graph over $\{1,2, \dots , n \}$ with no  matched edge (but the appropriate number of legs). For $1 \leq m \leq S_{ \mathrm{out}}$, we use the (matched) edges of $\mathrm{RM}( \mathbf{i}, \mathbf{a},m-1)$ to (try to) park the $m$th car which arrives on the vertex incident to the out-leg $X_m$. If we manage to park it, we denote by $\zeta_{m} \in \{1,2, \dots , n \}$ its parking spot, otherwise we set $\zeta_{m} = \dagger$. When $\zeta_{m} \ne \dagger$, this means that there is no edge going out from $\zeta_{m}$ but there is an unmatched out-leg since the underlying graph $\mathrm{RM}( \mathbf{i})$ is a mapping. 
We use the unmatched in-leg (corresponding in the mapping to) $ Y_m$ given in the construction of the frozen configuration model and match it with the out-leg going out from the white part of  $\zeta_{m}$ to form $\mathrm{RM}( \mathbf{i}^{(n)}, \mathbf{a}^{(n)},m)$.  We write this edge $\zeta_m \to Y_m$ (vertex $\to$ in-leg) .

\begin{center}
\includegraphics[width=12.5cm]{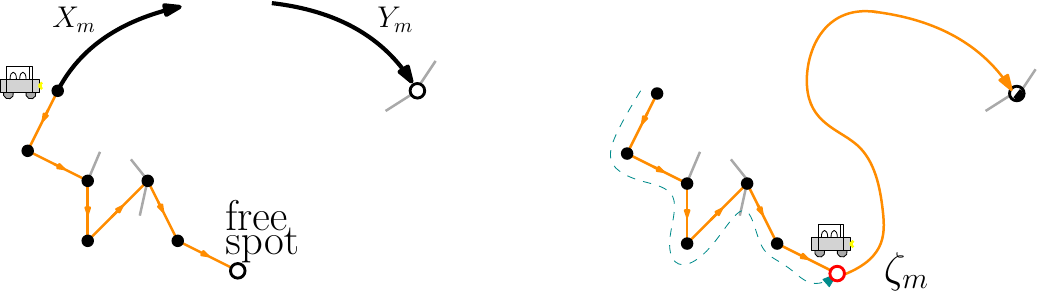}
\end{center}
The submapping $ \mathrm{RM}( \mathbf{i}^{(n)}, \mathbf{a}^{(n)}, S_{ \mathrm{out}})$ may have unmatched in-legs and vertices without out-edges remaining. Thus, in a second part, we match them uniformly at random to get $ \mathrm{RM}( \mathbf{i}^{(n)})$. See Figure \ref{fig:frozensteps} for a step-by-step illustration.
		\end{minipage}}
	
	\end{center}
	
\begin{figure}[!h]
 \begin{center}
 \includegraphics[width=15cm]{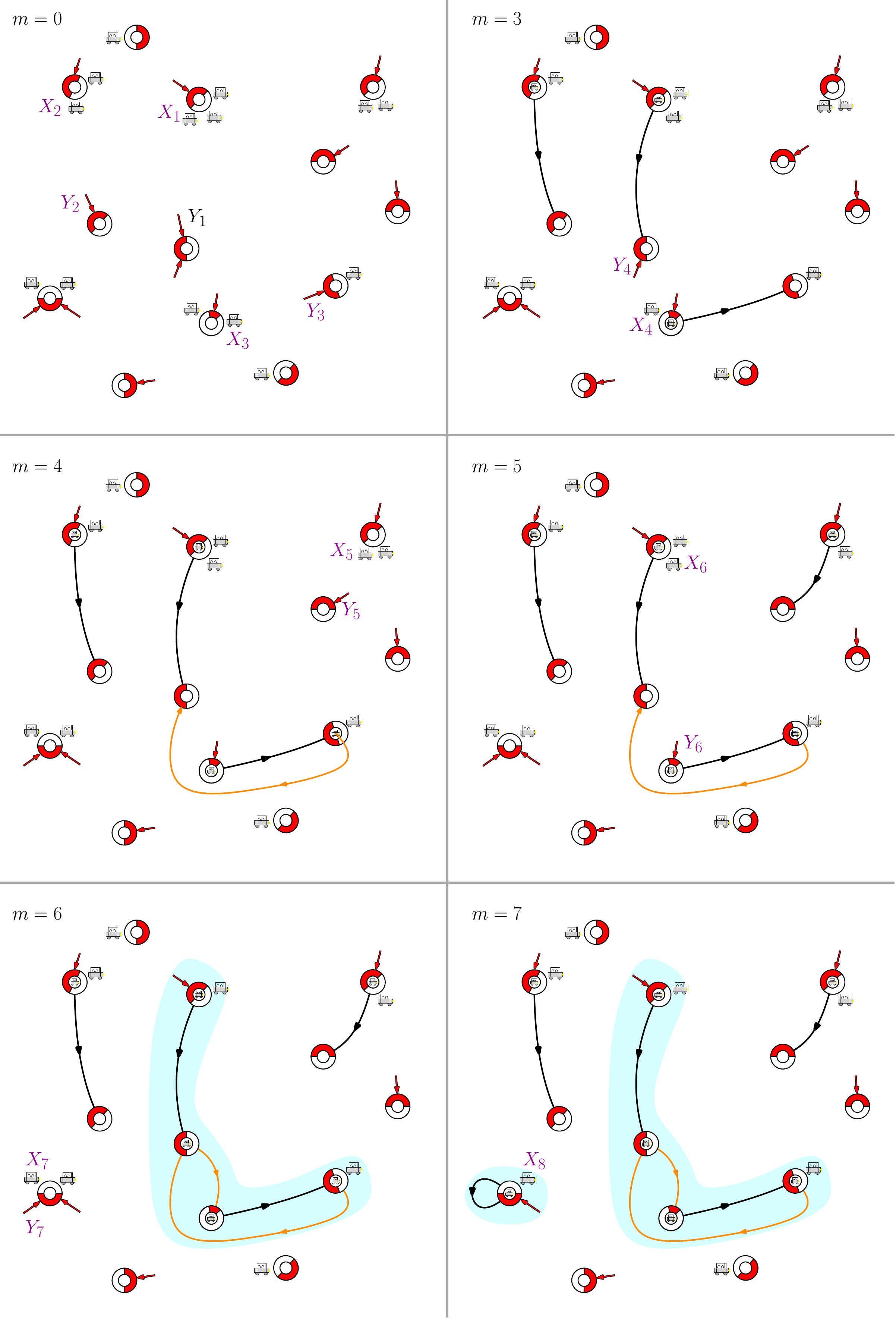}
 \end{center}
 \end{figure}
 \clearpage 
\newgeometry{a4paper,vmargin={1.5cm,2.5cm},hmargin={2cm,2cm}}
\begin{figure}[!h]
 \begin{center}
   \includegraphics[width=15cm]{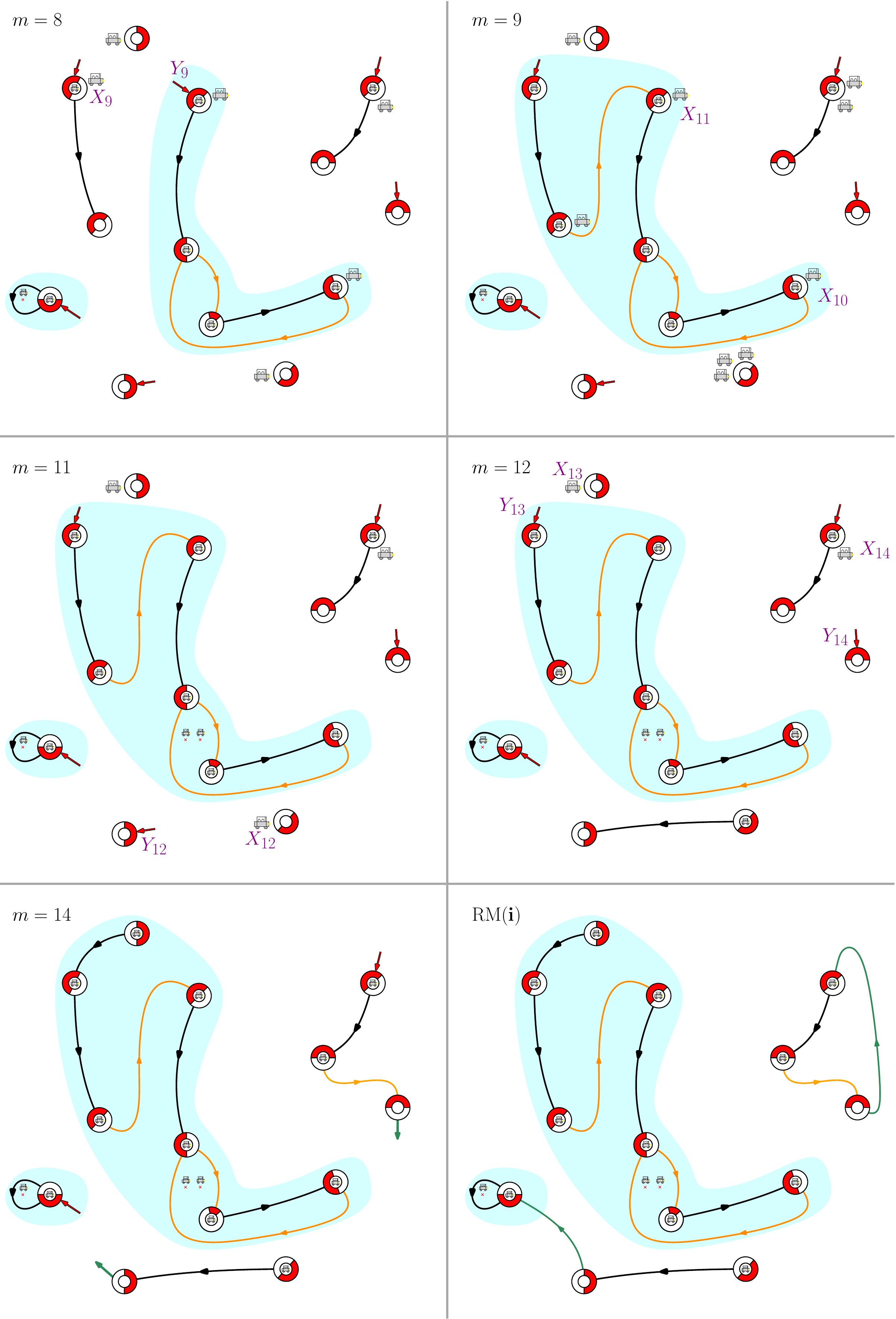} 
 \caption{ \label{fig:frozensteps}{Step-by-step illustration of the construction of the random mapping. In black, the edges that are also present in the construction of the frozen conguration model and in orange, the edges that are changed. In green in the end, the edges that come in the second step of the construction (not from the cars).}}
 \end{center}

 \end{figure}

\restoregeometry

\clearpage
\restoregeometry

It is crystal clear from the construction that the vertices having an outgoing edge in $ \mathrm{RM}( \mathbf{i}, \mathbf{a},m)$ are those which already accommodate a car at step $m$, and that if a car is caught in an endless loop, the submapping does not evolved. The main point is to check that the mapping induced by this construction has the appropriate law and that the car arrivals are independent of the constructed mapping. 
Observe that a way to obtain a mapping with the distribution of  $ \mathrm{RM}(\mathbf{i})$ is to fix a (deterministic) order of the vertices and one-by-one, for each vertex, choose an unmatched in-legs uniformly at random, which we match with an leg going out from the vertex defined by our (fixed) order.

Here, the construction follows almost the same lines with the small difference that the order of the vertices are random, but it is at each step independent of the already matched edges. Indeed, first note that conditionally on  $\mathrm{RM}( \mathbf{i}, \mathbf{a},m-1)$, the unmatched in-edge  $Y_m$ is uniform among the unmatched in-edges of $\mathrm{RM}( \mathbf{i}, \mathbf{a},m-1)$. The vertex $ \zeta_{m}$ is a mesurable function of   $ \mathrm{RM}( \mathbf{i}, \mathbf{a},m-1)$ and $X_m$, and is thus independent of $ Y_m$.  Thus, the first part of the algorithm gives an order of the some vertices, which is completed by the second part of the construction and at each step, they are assigned an independent random uniform target out-leg. Lastly, notice that the target edges are independent of the chosen order of the vertices.

\end{proof}

The statement of  Theorem \ref{thm:coupling} in the case of trees is similar to Proposition \ref{prop:couplingmapping}, but to prove it, we need to introduce a step-by-step construction of trees with prescribed degree sequence, which was straightforward in the case of mappings (Section \ref{sec:parkingmapping}).

\subsection{Peeling random trees with prescribed in-degree sequence} 
We present here the Markovian exploration of random trees with given degrees which is adapted from that of Cayley trees introduced in \cite{contat2022surprising,contat2021parking}. Recall that we consider plane trees with a prescribed number of children for each labeled vertex given by the sequence $ \mathbf{i}=(i_1, \dots, i_n)$. It is actually very similar to the case of mappings with prescribed in-degree sequence. The only difference is that in the case of trees, there exist a root without out-edge or parent. Once we fix an (deterministic) order of the vertices, a tree $ \mathbf{t}$ is uniquely determined by the sequence of different and distinguishable target-leg of each vertex, which is the edge going from a vertex to its parent, with the exception of the root in which case we say that its target-edge is $ \emptyset$. In what follows, we simply write \emph{target} instead of target-edge.

\begin{figure}[!h]
 \begin{center}
 \includegraphics[width=13cm]{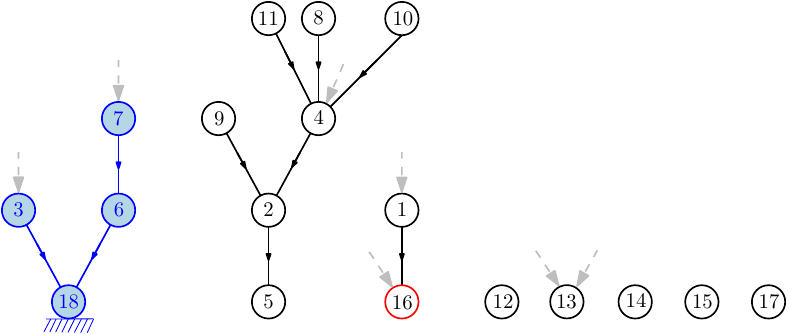}
 \caption{\label{fig:forest-ex}Exemple of a forest which we can observe during the peeling exploration of a tree}
 \end{center}
 \end{figure}

Thus, exploring a tree can be seen as revealing those $n$ targets one by one. Of course, not all sequences of targets correspond to a tree. Any subsequence of a sequence corresponding to a tree can be seen as a forest of rooted trees with some vertices which may contains unmatched legs, see Figure \ref{fig:forest-ex}. When the target of the root is “revealed” (one should probably better say that the root vertex is revealed) then we record this information by coloring the corresponding tree in blue, the other trees being referred to as white. 

Of course, a tree can be explored in $n!$ different ways depending on which of the $n!$ order of the vertices we choose. We shall  choose such an order using a  function $ \mathcal{A}$, called the \emph{peeling algorithm}, which to each forest $ \mathbf{f}$ (containing blue tree or not), associates a vertex $ \mathcal{A}( \mathbf{f})$ whose target is not yet revealed. The peeling of a tree $ \mathbf{t}$ with algorithm $ \mathcal{A}$ is then the  sequence 
$$ (\mathbf{f}_0^ \mathcal{A}, \mathbf{f}_1^ \mathcal{A}, \dots , \mathbf{f}_n^ \mathcal{A} = \mathbf{t} )$$ 
where $\mathbf{f}_0^ \mathcal{A}$ is the forest with $n$ isolated vertices with the unmatched legs prescribed by $ \mathbf{i}$ and of all $m \leq n-1$, the forest $ \mathbf{f}_{m+1}^ \mathcal{A}$ is obtain from the forest $ \mathbf{f}_m^ \mathcal{A}$ by adding an edge from the vertex $ \mathcal{A}( \mathbf{f}_m^{ \mathcal{A}})$ to its target-leg present in $ \mathbf{t}$. This exploration is called a peeling exploration and as in the case of Cayley trees, when the underlying tree is a uniform random tree with prescribed degree, this exploration is a Markov chain with explicit probability transitions. The following proposition is the analogue of \cite[Proposition 6]{contat2021parking}.

\begin{proposition}[Markov]\label{prop:peel} Fix a peeling algorithm $ \mathcal{A}$ and an in-degree sequence $ \mathbf{i}= (i_1, \dots, i_n)$. If $T( \mathbf{i})$ is a uniform random tree with degree given by $ \mathbf{i}$, then the peeling exploration $(\mathbf{f}_m^ \mathcal{A} : 0 \leq m \leq n)$ of $T( \mathbf{i})$ with peeling algorithm $ \mathcal{A}$ is a Markov chain whose transitions are described as follows. Conditionally on  $\mathbf{f}_m^ \mathcal{A}$ and on $ \mathcal{A}(\mathbf{f}_m^ \mathcal{A})$, we denote by $k \geq 0$ the number of unmatched edges in the tree with root $ \mathcal{A} ( \mathbf{f}_m^ \mathcal{A})$ and by $ \ell \geq 0$ the number of unmatched edges in the blue tree (if any). Then, 
\begin{itemize}
\item If $ \ell = 0$, with probability $ \frac{k}{n-1-m}$, then $ \mathcal{A}( \mathbf{f}_m^ \mathcal{A})$ is the root of the underlying tree $ T( \mathbf{i})$ and we color the tree of $\mathbf{f}_m^ \mathcal{A}$ in blue; otherwise, the target-edge of $ \mathcal{A} (\mathbf{f}_m^ \mathcal{A})$ is a uniform edge not belonging to the tree of root $ \mathcal{A} (\mathbf{f}_m^ \mathcal{A})$.
\item if $ \ell \geq 1$, with probability $ \frac{ \ell-1 +k}{n-1-m}$, the target edge of $ \mathcal{A} ( \mathbf{f}_m^ \mathcal{A})$ is a uniform unmatched edge of the blue tree of $ \mathbf{f}_m^ \mathcal{A}$, otherwise it is a uniform vertex of the remaining tree except the tree of root $ \mathcal{A}( \mathbf{f}_m^ \mathcal{A})$.
\end{itemize}
\end{proposition}

Note $ \ell \geq 1$ if and only if the root is already revealed but not all edges of the tree ($ m \neq n$), and that if $m<n$, the number of unmatched edges is $n-1-m$ if $ \ell =0$, and $n-m$ otherwise. As in \cite{contat2021parking} and \cite{contat2022surprising}, the proof relies on counting formulas, see e.g. \cite{harary1964number}. 

\begin{lemma} If $ \mathbf{f}$ is a forest of white rooted trees and prescribed degrees with $ k$ trees, then the number of rooted trees containing $ \mathbf{f}$ is $(k-1)!$. If $ \mathbf{f}^{ \star}$ is a forest of white rooted trees and prescribed degrees with $ k$ trees including a blue tree with $ \ell \geq 1$ unmatched edges, then the number of rooted trees containing $ \mathbf{f}^ \star$ is $ \ell (k-2)!$.
\end{lemma}

Equipped with this lemma, we are now able to prove Proposition \ref{prop:peel}.
\proof[Proof of Proposition \ref{prop:peel}] The proof is a direct consequence of the previous lemma. It suffices to notice that  since the underlying tree $ T ( \mathbf{i})$ is uniform over all tree with degrees given by $ \mathbf{i}$, then for all $m \geq 0$, conditionally on $\mathbf{f}_m^ \mathcal{A}$, the tree  $ T ( \mathbf{i})$ is a uniform tree among those which contain the forest  $\mathbf{f}_m^ \mathcal{A}$ (with or without a blue tree depending whether the root vertex has been revealed or not). Hence, for every (compatible) target $e$ (possibly $ \varnothing$), the probability that $\mathbf{f}_{m+1}^ \mathcal{A}$ is obtained from $\mathbf{f}_m^ \mathcal{A}$ by adding an edge from $ \mathcal{A}(\mathbf{f}_m^ \mathcal{A})$ to $e$ (or revealing that $ \mathcal{A} (\mathbf{f}_m^ \mathcal{A})$ is the root if $e= \varnothing$) is 

$$ \frac{ \# \{ \mathbf{t} \mbox{ containing }\mathbf{f}_m^ \mathcal{A}\mbox{ and }\mathcal{A} (\mathbf{f}_m^ \mathcal{A})  \to e \} }{\# \{ \mathbf{t} \mbox{ containing }\mathbf{f}_m^ \mathcal{A} \} }.$$
Using the previous lemma, we recognize the transition probabilities given in Proposition \ref{prop:peel} and obtain the desired result. \endproof

Another step-by-step construction of trees with prescribed out-degree sequence is present in \cite{addario2023foata}. However, the authors consider trees where the leaves are labeled independently of the other vertices. 

\subsection{Coupling between parking on trees and the frozen configuration model}\label{sec:couplingtree}

Equipped with this peeling exploration of trees with a prescribed degree sequence, we have now all the tools to adapt our coupling from the case of mappings to the case of trees. The main difference with the case of mappings is that we will need to add an additional source of randomness. 

\begin{proposition}\label{prop:couplingtree} Let us fix $ n \geq 1$, an in-degree sequence $ \mathbf{i}$ such that $ S_ \mathrm{in} := \sum_{k \geq 1} i_k =n-1$ and an out-degree sequence $ \mathbf{a}$. 
We can couple the parking process on $ T ( \mathbf{i})$  with car arrivals given by $ \mathbf{a}$ and the frozen configuration model $ \mathrm{FCM} ( \mathbf{i}, \mathbf{a})$ such that 
\begin{itemize}

\item If the root contains a car in $T (\mathbf{i})$, then its parked component is made of all blue vertices in $ \mathrm{FCM} ( \mathbf{i}, \mathbf{a})$. All other near parked components in $T ( \mathbf{i})$  are the white connected components of $ \mathrm{FCM} (\mathbf{i},\mathbf{a})$ (in terms of subset of vertices).

\item the out-legs deleted while constructing $ \mathrm{FCM} ( \mathbf{i}, \mathbf{a})$ correspond to the unparked cars in $ T (\mathbf{i})$.
\item  the tree $ T ( \mathbf{i})$ is a uniform tree among those containing the forest of the nearly parked components.
\end{itemize}

\end{proposition}

\proof As in the case of mappings, we construct the tree $ T ( \mathbf{i})$ using the edges of the frozen configuration model. The main difference is that the appearance of the first cycle in the construction of $ \mathrm{FCM}( \mathbf{i}, \mathbf{a})$ corresponds to the detection of the root in the exploration of the tree, and we then need an additional source of randomness to redirect some edges to keep the tree structure.

\begin{center}
	\fbox{\begin{minipage}{15cm}
			\paragraph{From oriented edges to parking on trees.} 
Once again, the out-legs are interpreted as the car arrivals on the tree $ T( \mathbf{i})$. We use them and their order given by the sequence $(X_m : 1 \leq m \leq S_{ \mathrm{out}})$  to construct iteratively an ``growing" sequence of forests $ ( F_{m} ^{ \mathrm{park}} : 0 \leq m \leq S_ \mathrm{out})$ with zero or one blue tree. Note that this time, we have to be more careful when mentioning the in-legs since they will not exactly all have the same status (matched or unmatched) in the frozen configuration model and in the forest $ \mathbf{F}_{ m }^{ \mathrm{park}}$.

Initially $ F_{0}^{\mathrm{park}}$ is the forest with $n$ isolated vertices with the unmatched in-legs prescribed by $ \mathbf{i}$. For $m \geq 1$, we use the edges of $ F_{m-1} ^{\mathrm{park}}$ to (try to) park the $m$th car which arrive on the vertex incident to the leg $X_m$. If we manage to park it, we denote by $\zeta_{m} \in \{1,2, \dots , n \}$ its parking spot, otherwise set $\zeta_{m} = \dagger$. If $\zeta_{m}= \dagger$ then $ F_{m}^{ \mathrm{park}} = F_{m-1}^{ \mathrm{park}}$. Otherwise, we denote by $Y_m$ the unmatched in-leg given at step $m$ in the construction of the frozen configuration model $\mathrm{FCM}( \mathbf{i}, \mathbf{a})$. Then, 

	\begin{itemize}
			\item if $Y_m$ is incident to a white vertex and the addition of the edge $X_{m} \to Y_{m}$ does not create a cycle in the frozen configuration model construction, then $Y_m$ is also unmatched in $ \mathbf{F}_{m-1}^{ \mathrm{park}}$ and we add the edge $\zeta_{m} \to Y_{m}$ to $ \mathbf{F}_{m-1}^{ \mathrm{park}}$ to form $ \mathbf{F}_{ m }^{ \mathrm{park}}$,
			\item if the addition of the edge $X_{m} \to Y_{m}$ creates a cycle in the frozen configuration model construction or if the vertex incident to $Y_m$ is blue, then 
			\begin{itemize}
			\item If $F_{m-1}^{ \mathrm{park}}$ has no blue tree (the root vertex is not revealed), then color the tree of $ \zeta_m$ in blue to form $  F_{m} ^{\mathrm{park}}$, to reveal that $ \zeta_m$ is the root of the tree. The leg $Y_m$ is not matched in $  F_{m} ^{\mathrm{park}}$ but it is matched in the frozen configuration model. Thus we can note re-use it as another $Y_k$ for $k> m$, but it can be used after to construct the tree $T( \mathbf{i})$.
			\item Otherwise add  $\zeta_{m} \to U_{m}$ where $U_{m}$ is a uniform unmatched in-leg over the blue tree of $  F_{m-1} ^{\mathrm{park}}$ sampled independently of the past to form $  F_{m} ^{\mathrm{park}}$.
\end{itemize}
			\end{itemize}

		\end{minipage}}
	\end{center}

As in the case of mapping, the forest $ \mathbf{F}_{ S_{ \mathrm{out}} }^{ \mathrm{park}}$ at step $S_{ \mathrm{out}}$ is a subforest of the underlying tree $ T ( \mathbf{i})$, which corresponds to its near parked components, and the deterministic properties of the coupling between the construction of $ \mathrm{FCM}( \mathbf{i}, \mathbf{a})$ and the parking process on $ T ( \mathbf{i})$ work exactly in the same way. 

It only remains to check that, if we choose a (deterministic) order of the roots of the white trees of  $ \mathbf{F}_{ S_{ \mathrm{out}} }^{ \mathrm{park}}$ and sample their target-edges using the transition described in Proposition \ref{prop:peel}, we obtain a uniform tree with prescribed degrees, independent of the sequence $ \mathbf{a}$.  To see this, we interpret the Markov chain $( \mathbf{F}_{m}^{ \mathrm{park}} : m \geq 0)$ as the first steps of an exploration of a uniform tree $ T ( \mathbf{i})$ with numbers of children prescribed by $ \mathbf{i}$. We construct the peeling algorithm $ a_{ \mathrm{park}}$ as follows. At $m=0$, we start from the forest with isolated vertices, and for $m \geq 1$, we let a car arrive on vertex $X_m$ and the car follows the edges present in $ F_{m-1}^{ \mathrm{park}}$ to park in $ \zeta_m$. If the car does not park, then we put $ \zeta_m = \dagger$ and we move to step $m+1$. If $ \zeta \neq \dagger$, we then put $$ a_{ \mathrm{near}} (F_{m-1}^{ \mathrm{park}}) = \zeta_m, $$
and we reveal the target-edge $Y_m$ of $ \zeta_m$, that is, we add the edge $ \zeta_m \to Y_m$ to form $ F_m^{ \mathrm{park}}$. Then the process $( \mathbf{F}_{m}^{ \mathrm{park}} : m \geq 0)$ as the law of the exploration $(F_m^{ \mathrm{park}} : m \geq 0)$. Moreover, the peeling algorithm $a_{ \mathrm{near}}$ is a measurable function of the $X_m$'s as well as the second part of the peeling exploration. Thus by Proposition \ref{prop:peel}, the tree $ T( \mathbf{i})$ constructed this  way is a uniform tree with prescribed number of children, and in particular, the tree is independent of $ \mathbf{a}$. The proposition follows.

\endproof

\section{Phase transition via the local limit}

In this section, we introduce even more randomness and suppose that the given degrees and car arrivals are random. We suppose that we work under the hypothesis \eqref{eq:hypo-convergence} with $ \E_{{\nu}} [m] \leq 1$. Recall from this hypothesis the notation for the sequence of in-degree sequences $ \mathbf{I}^{(n)}$.

\subsection{Local limit of trees and mappings with given degree sequences}
In their proof of the location of the phase transition for parking on Bienaym\'e--Galton--Watson trees, Curien and Hénard \cite{curien2022phase} used a differential equation on the expectation of the flux, which is obtained via a spine decomposition and the local limit of these trees to characterize the critical point. Our proof of Theorem \ref{thm:coupling:curienhenard} is more ``geometric" since it relies on the apparition of a giant component in the frozen configuration model or on the configuration model and a soft  continuity argument for the probability to observe a macroscopic flux of outgoing cars with respect to the local topology. We thus analyse the local limit of the components of parked cars using our coupling with random graphs.

\paragraph{Local limits and Benjamini-Schramm quenched convergence.} Let us start by reminding the notion of local convergence that we require. Recall that we allow our graphs to have unmatched (oriented) legs.

\begin{definition} Let $(G_n)$ be a sequence of finite (random) graphs and for all $n \geq 1$, conditionally on $ G_n$, we sample $X_n$ and $Y_n$ two vertices of $G_n$ chosen independently and uniformly at random, and $ (G_{ \infty}, \rho_{ \infty})$ a (random) graph rooted at $\rho_{ \infty}$. 
We say that $(G_n)$ converges Benjamini-Schramm quenched towards $(G_ \infty, \rho_ \infty)$ if for all $ r \geq 0$ and for all non negative functions $f$ and $g$ such that the value of $ f ( \mathfrak{g}, \rho)$ and $ g ( \mathfrak{g}, \rho)$ only depends on the ball $ B_r ( \mathfrak{g}, \rho)$ of $ \mathfrak{g}$ of radius $ r$ around $ \rho$, 
$$ \mathbb{E} \left[ f (G_n, X_n) g (G_n, Y_n) \right ] \xrightarrow[n\to\infty]{} \mathbb{E} \left[ f (G_ \infty, \rho_ \infty) \right] \mathbb{E} \left[ g (G_ \infty, \rho_ \infty) \right ].$$
\end{definition}

\begin{figure}[!h]
 \begin{center}
 \includegraphics[width=14cm]{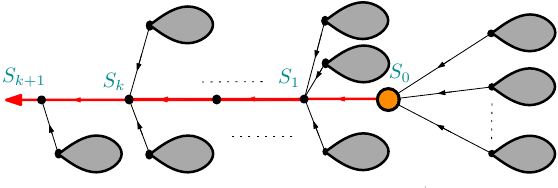}
 \caption{\label{fig:sintree}Illustration of the construction of the sin-tree of Aldous. In orange, the root vertex. Each gray tree is a Bienaym\'e--Galton--Watson tree with offspring distribution $ \nu$. The number of children of the vertex $S_0$ has law $ \nu$ whereas the number of children of the vertices $S_k$ with $ k \geq 1$ has law $ \overline{ \nu}$.}
 \end{center}
 \end{figure}

Recall that in our case, we are under hypothesis \eqref{eq:hypo-convergence} and in particular, the empirical distribution of the in-degree sequence converges towards $ \nu$. Heuristically, the tree $ T ( \mathbf{I}^{(n)})$ is not far from a Bienaym\'e--Galton--Watson tree with law $ \nu$ conditioned to have size $n$. Thus its Benjamini-Schramm quenched limit should be the sin-tree of Aldous \cite{aldous1991asymptotic} with law $\nu$, which we describe now. We consider a semi-infinite line $S_0, S_1, \ldots$ called the spine,  which we see rooted at $ \rho_ \infty = S_0$, and for $ i \geq 1$, we orient the edges from $S_{i-1}$ to $S_{i}$ and graft independently on each $S_i$ a random number number $Y-1$ of independent Bienaym\'e--Galton--Watson trees with offspring distribution $ \nu$ where $ Y \sim \overline{\nu}$, and consider a random uniform ordering of the children of $ S_i$ and furthermore, graft independently $ X \sim \nu$ Bienaym\'e--Galton--Watson trees on $S_0$. We obtain a tree rooted which we denote by $ (T_{ \infty} (\nu), \rho_ \infty)$,  see Figure \ref{fig:sintree}.

\begin{proposition}\label{prop:limloc} Under the hypothesis \eqref{eq:hypo-convergence},  the random trees sequence $ (T ( \mathbf{I}^{(n)}))$ and the random mappings sequence $ (\mathrm{RM}( \mathbf{I}^{(n)}))$ converge Benjamini--Schramm quenched towards the sin-tree  $ (T_{ \infty} (\nu), \rho_ \infty)$.
\end{proposition}

\proof First, thanks to \cite[Proposition 11]{aldous1991asymptotic}, it suffices to consider the ``top-part" of the trees and mapping rooted at the uniform point. More precisely, let us focus on the case of the tree $ T ( \mathbf{I}^{(n)})$ and let $X_n$ and $Y_n$ be two uniform and independent points of  $ T ( \mathbf{I}^{(n)})$. Given an oriented graph $ \mathfrak{g}$, we denote by $ \mathrm{Top}( \mathfrak{g}, x)$ the (plane) graph of the descendants of $x$ in $ \mathfrak{g}$ (without labelling of the vertices) i.e.\ the graph made of all vertices $v$ such that there exists a \emph{direct path} from $v$ to $x$. We only need to prove that for all fixed rooted tree $ \mathfrak{t}_1$ and $ \mathfrak{t}_2$, we have $$ \mathbb{P} \left( \mathrm{Top } (T ( \mathbf{I}^{(n)}), X_n) = \mathfrak{t}_1 \mbox{ and } \mathrm{Top } (T ( \mathbf{I}^{(n)}), Y_n) = \mathfrak{t}_2\right) \xrightarrow[n\to\infty]{} \mathbb{P} ( \mathcal{T} = \mathfrak{t}_1) \mathbb{P} ( \mathcal{T} = \mathfrak{t}_2),$$ 
where $ \mathcal{T}$ is a Bienaym\'e--Galton--Watson tree with offspring distribution $ \nu$. This is in fact quite straightforward and let us start for the sake of clarity with the case of one fixed tree $ \mathfrak{t}_1$. We denote by $k_i$ the number of vertices of degree $i$ in $ \mathfrak{t}_1$, and $k$ its total number of vertices.  To compute the probability that $ \mathrm{Top } (T ( \mathbf{I}^{(n)}), X_n) = \mathfrak{t}_1 $, we sum over all possible labellings of the vertices of $ \mathfrak{t}_1$ and compute the probability of the connection of the edges, which is simply $ \prod_{i = 1}^{k-1} 1/(n-i)$. Combining this with the probability that $X_n$ has the good label, we obtain 
\begin{align*}  \mathbb{P} (  \mathrm{Top } (T ( \mathbf{I}^{(n)}), X_n) = \mathfrak{t}_1) &= \sum_{\mbox{possible labelings}} \prod_{i = 0}^{k-1} \frac{1}{(n-i)} \\
&= \left( \frac{1}{n} + o \left( \frac{1}{n}\right)\right)^k \cdot \prod_{i \geq 0} \binom{\# \{j\ \mbox{s.t.}\ i_j = i\}}{k_i} \\
&= \left( \frac{1}{n} + o \left( \frac{1}{n}\right)\right)^k \cdot \prod_{i \geq 0} \left( n \nu_i + o(n)\right)^{k_i} \\
& \xrightarrow[n\to\infty]{}\prod_{i \geq 0} \nu_i^{k_i} = \mathbb{P} ( \mathcal{T} = \mathfrak{t}_1),  \end{align*}
since $\nu_i$ is the asymptotic proportion of vertices of degree $i$ in $T( \mathbf{I}^{(n)})$. This computation can be easily extended to the case of two points $X_n$ and $Y_n$ and two trees $ \mathfrak{t_1}$ and $ \mathfrak{t}_2$, and we obtain the same probabilities in the case of random mappings. In particular, the descendants of a uniform point in a mapping form a tree in the limit, and a uniform point is not in an oriented cycle with high probability. We left this to the reader and conclude the proof.  \endproof

As a corollary, we can adapt Theorem 4.1 in \cite{contat2023parking} to obtain a phase transition result for the parking process on $  (T ( \mathbf{I}^{(n)}))$ and  $ (\mathrm{RM}( \mathbf{I}^{(n)}))$, which only depends on the Benjamini--Schramm limit.

\begin{corollary}\label{cor:phasetransition}We suppose \eqref{eq:hypo-convergence} and that there exists a constant  $K$ such that $m_{(k)}< K$ and $ \sigma_{(k)}^2< K$ for all  $ k \geq 0$. We then have 
\begin{eqnarray*} \frac{ \varphi ( T ( \mathbf{I}^{(n)}))}{n} \underset{ n \to \infty} {\xrightarrow{( \mathbb{P})}} C \qquad \& \qquad  \frac{ \varphi ( \mathrm{RM}( \mathbf{I}^{(n)}))}{n} \underset{ n \to \infty} {\xrightarrow{( \mathbb{P})}}
C 
\end{eqnarray*} where $C>0$ if and only if $ \mathbb{P}\left( \mbox{there exists an infinite cluster of parked cars in } \mathcal{T}_ \infty (\nu)\right) >0$.
\end{corollary}

Actually, Theorem 4.1 in \cite{contat2023parking} only consider a sequence of \emph{trees} which converges Benjamini--Schramm quenched but there is no obstacle in the proof to extend the result to a sequence of random mappings converging Benjamini--Schramm quenched towards a tree (with one end).

\subsection{Local limit of the frozen configuration model}
\begin{figure}[!h]
 \begin{center}
 \includegraphics[width=10cm]{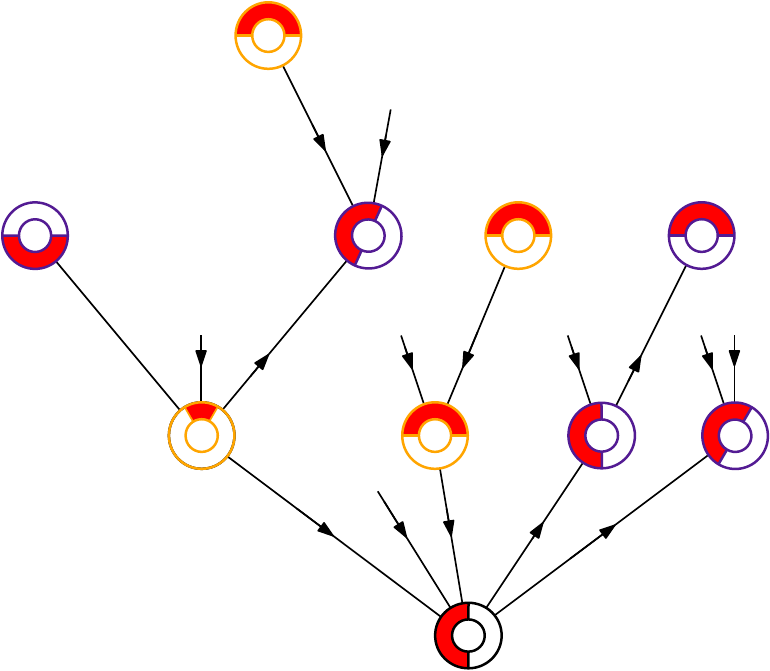}
 \caption{\label{fig:limlocpark}Illustration of a possible realization of the generalized multi-type Bienaym\'e--Galton--Watson tree. All vertices are white for the frozen configuration model but we used colors to distinguish their type in the Bienaym\'e--Galton--Watson tree: the root vertex is in dark, the vertices of type $o$ are in orange and the vertices of type $i$ are in indigo. Each vertex have two sides: red for in in-legs and white for the out-legs.}
 \end{center}
 \end{figure}

The goal of this section is to investigate the local limit of the frozen configuration model, so that using Corollary \ref{cor:phasetransition}, we can characterize the phase transition of the parking process. First let us describe the multi-type Bienaym\'e--Galton--Watson trees which are involved in the local limit of the oriented configuration model where the degree distribution is given by $\lambda$ in \eqref{eq:hypo-convergence}. Recall that we explore the weakly connected component of the graph starting from a point uniformly at random. As for many sparse random graph model, the connected component of a uniform point is locally a tree. In our case, the tree has three types of vertices : the root vertex has its own type, and the other vertices are either of type $ o$ (with a size-biased number of out-edges) or of type $ i$ (with a size-biased number of in-edges), see Figure \ref{fig:limlocpark}. Recall that each vertex has two sides, a white one for the out-legs and a red one for the in-legs.

\begin{definition}
Suppose that $ \mathbb{E}_{\nu} [m] \leq 1$. Let us consider the multi-type Bienaym\'e--Galton--Watson tree where all vertices reproduce independently with the following offspring distributions: 

\begin{itemize}
\item We start for a root vertex $ \varnothing$ which has its own ``normal type". We sample a variable $X_{ \varnothing}$ according to $ \nu$, which is the typical number of in-legs of a vertex and conditionally on $X_{ \varnothing}$, we sample a variable $ Y_{ \varnothing}$ according to $ \mu_{ (X_{ \varnothing})}$, which is the typical number of out-legs of this vertex. However, not all in-legs are matched in the end, thus, for $1 \leq k \leq X_{ \varnothing}$, we either add  a child of type $o$ (chosen  with a size biased number of out-legs) with probability $ \mathbb{E}_{ \nu} [m]$, or an unmatched in-leg. We also add $Y_{ \varnothing}$ children of type $i$.
\item A vertex $v$ of type $o$ has a size-biased number of out-legs. Thus, we sample a pair of variables $(X_v, Y_v)$ according to 
$$ (X_v, Y_v) \sim \frac{1}{ \mathbb{E}_{ \nu} [m]}\sum_{k \geq 0} \sum_{ j \geq 0} j \nu_k \mu_{(k),j} \delta_{(k,j)}. $$
The second variable $Y_v$ gives the total number of out-legs and thus $v$ has $ Y_v - 1$  children of type $i$. We choose a corner uniformly at random between two out-legs (there is $Y_v$ possibilities) and put the red side of the vertex in this corner. In this corner from left to right, for $1 \leq k \leq X_v$, we either add a child of type $o$ with probability $ \mathbb{E}_{ \nu} [m]$, or an unmatched in-leg $ 1- \mathbb{E}_{ \nu} [m]$.
\item A vertex of type $i$ has a size-biased number of in-legs. Thus, we sample a pair of variables $(X_v,Y_v)$ according to 
$$ (\overline{ \nu}, \mu_{ \overline{ \nu}}) := \sum_{k \geq 0} \sum_{j \geq 0} k \nu_k  \mu_{(k),j} \delta_{(k, j)}.$$  
 For $1 \leq k \leq X_v -1$, we either add a child of type $o$ with probability $ \mathbb{E}_{ \nu} [m]$, or an unmatched in-leg with probability $ 1- \mathbb{E}_{ \nu} [m]$. We then choose a corner uniformly at random between two in-legs (there is $X_v$ possibilities) and put the white side of the vertex in this corner where we add $Y_v$ out-legs. 
\end{itemize}
We denote by $ \mathcal{T}$ a random tree (with possibly remaining unmatched legs) with this law.
\end{definition}

 First note that  when  $\E_{\overline{\nu}} [m] \leq 1$ and $ \Theta \geq 0$, this tree $ \mathcal{T}$ is almost surely finite. Indeed, let us compute the expected number of children of each type $i$ and $o$ for each type of vertices. For $v_o$ a vertex of type $o$, its expected number of children of type $o$ is 
 \begin{eqnarray*}\mathbb{E} [ Y_{v_o} - 1]  &=&  \frac{1}{ \mathbb{E}_{ \nu} [m]} \left(\sum_{k \geq 0} \sum_{ j \geq 0} j^2 \nu_k \mu_{(k),j} \right) - 1 \\
 &=&\frac{1}{ \mathbb{E}_{ \nu} [m]} \left(\sum_{ k \geq 0}  \nu_k \left(\sigma^2_{(k)} + m^2_{(k)} \right) -  \mathbb{E}_{ \nu} [m] \right) = \frac{ \mathbb{E}_{ \nu}[ \sigma^2 + m^2 - m] }{ \mathbb{E}_{ \nu} [m]}.
 \end{eqnarray*}
Similarly, its expected number of type $i$ is 
 \begin{eqnarray*}\mathbb{E} [ X_{v_o} ] \cdot   \mathbb{E}_{ \nu} [m] =  \frac{ \mathbb{E}_{ \nu} [m] }{ \mathbb{E}_{ \nu} [m]} \sum_{k \geq 0} \sum_{ j \geq 0} k j  \nu_k \mu_{(k),j} 
 = \sum_{ k \geq 0} k \nu_k  m_{(k)}=  \mathbb{E}_{  \overline{\nu}}[  m].
 \end{eqnarray*}
With very similar computations for the expected number of children of a vertex of type $i$, we obtain the following mean matrix
\begin{equation*} 
\begin{pmatrix} 
\mathbb{E}_{ \overline{ \nu}} [ m ]  & \Sigma^2 \mathbb{E}_{ \nu}[m]   \\
\frac{ \mathbb{E}_{ \nu}[ \sigma^2 + m^2 - m] }{ \mathbb{E}_{ \nu} [m]}&\mathbb{E}_{ \overline{ \nu}}[m] 
\end{pmatrix},
 \end{equation*}
which has two eigenvalues $ \lambda_{\pm} = \mathbb{E}_{ \overline{ \nu}}[m] \pm \sqrt{  \Sigma^2 \mathbb{E}_{ \nu}[ \sigma^2 + m^2 - m]}$. Thus $\E_{\overline{\nu}} [m] \leq 1$ and $ \Theta \geq 0$ if and only if $ \lambda_{ \pm} \leq 1$ and in this case the Bienaym\'e--Galton--Watson tree $ \mathcal{T}$ is almost surely finite. 

Recall that we assume $\E_{{\nu}} [m] \leq 1$. 
We can show that in this regime, i.e.\ when $\E_{\overline{\nu}} [m] \leq 1$ and $ \Theta \leq 0$,   this tree is also the Benjamini--Schramm limit of the frozen configuration model. Note that this tree is infinite when $\E_{\overline{\nu}} [m] \geq 1$ or $ \Theta < 0$, and this regime corresponds to the apparition of a giant component for the configuration model. However this tree will not be the local limit of the frozen configuration model since more out-legs will be deleted.

\begin{proposition}\label{prop:souscrit}Suppose that  $\E_{\overline{\nu}} [m] \leq 1$ and $ \Theta \geq 0$. Under \eqref{eq:hypo-convergence} with $\E_{{\nu}} [m] \leq 1$, the frozen configuration model converges  Benjamini--Schramm quenched towards the Bienaym\'e--Galton--Watson tree $ \mathcal{T}$ described above which is almost surely finite.
\end{proposition}
In particular, when  $ \mathbb{E}_{\nu} [m] \leq 1$ and $ \Theta \geq 0$,  the Benjamini--Schramm quenched limit of $\mathrm{FCM}(  \mathbf{I}^{(n)}, \mathbf{A}^{(n)})$ is almost surely a finite white tree. As a consequence, the near parked component of the root in $ \mathcal{T}_ \infty ( \nu ) $ is almost surely finite and the parking process is subcritical by our coupling and Corollary \ref{cor:phasetransition}. 

\begin{proof}

We want to prove that the random tree $ \mathcal{T}$ is the appropriate Benjamini--Schramm quenched limit of $ \mathrm{FCM}( \mathbf{I}^{(n)}, \mathbf{A}^{(n)})$ when $ \Theta \geq 0$. First, when we look at the graph $ \mathrm{FCM}( \mathbf{I}^{(n)}, \mathbf{A}^{(n)})$ seen from a uniform vertex, this vertex is the ``root" of our graph and we put its in-legs or in-edges on the left and its out-edges on the right (see Figure \ref{fig:limlocpark}).
As in the proof of Proposition \ref{prop:limloc}, we only verify the convergence of the distribution of $ \mathrm{FCM}( \mathbf{I}^{(n)}, \mathbf{A}^{(n)})$ seen from one uniform point, but the extension to two points is straightforward. We thus fix a finite tree $ \mathfrak{t}$ with a root, with $k^{ \mathrm{out}}$ vertices of type $o$ (thus $k^{ \mathrm{out}}$ matched in-legs in total), with $A$ un-matched in-legs in total , and where $(d^{ \mathrm{in}}_v : v \in \mathfrak{t})$ and $(d^{ \mathrm{out}}_v : v \in \mathfrak{t})$ are respectively the  in-degree and out-degree sequences. We want to compute the probability that $ \mathcal{T} = \mathfrak{t}$. For each vertex $v$ of type $o$, the probability that it has the right degrees (with our Bienaym\'e--Galton--Watson distribution)  is $d^{ \mathrm{out}}_{v} \nu_{d^{ \mathrm{in}}_{v}} \mu_{(d^{ \mathrm{in}}_{v} ),d^{ \mathrm{out}}_{v}}/ \mathbb{E}_{ \nu} [m]$ and the probability to put the red part of the vertex at the right location between the out-edges is $ 1/d^{ \mathrm{out}}_{v}$, which makes $\nu_{d^{ \mathrm{in}}_{v}} \mu_{(d^{ \mathrm{in}}_{v} ),d^{ \mathrm{out}}_{v}}/\mathbb{E}_{ \nu} [m]$ in total. Similarly, for a vertex $v$ of type $i$, the probability the observe the right degree sequence with the good place for the white part is $\nu_{d^{ \mathrm{in}}_{v}} \mu_{(d^{ \mathrm{in}}_{v} ),d^{ \mathrm{out}}_{v}}$. We then obtain
\begin{equation*}  
\mathbb{P} \left( \mathcal{T} = \mathfrak{t}\right)= \frac{1}{\mathbb{E}_{ \nu} [m]^{k_o}}\left(\prod_{ v \in \mathbf{t}}  \nu_{d^{ \mathrm{in}}_{v}} \mu_{(d^{ \mathrm{in}}_{v} ),d^{ \mathrm{out}}_{v}} \right)\cdot \mathbb{E}_{ \nu} [m]^{k_o} (1- \mathbb{E}_{ \nu} [m])^A
= \left(\prod_{ v \in \mathbf{t}} \nu_{d^{ \mathrm{in}}_{v}} \mu_{(d^{ \mathrm{in}}_{v} ),d^{ \mathrm{out}}_{v}} \right)\cdot (1- \mathbb{E}_{ \nu} [m])^A.\end{equation*}

Now let us compute the probability to observe this tree $ \mathfrak{t}$ seen from a uniform vertex in the frozen configuration model. Let us call $k$ the total number of vertices of $\mathfrak{t}$ and write $(v_1, \dots, v_k)$ for the vertices of $ \mathfrak{t}$ in the depth-first order. Again, we sum over all possible labelings of the vertices of $ \mathfrak{t}$ and compute the probability of the edge-connection. Thus conditionally on $ \mathbf{I}^{(n)}, \mathbf{A}^{(n)}$, the probability that $ \mathrm{FCM}( \mathbf{I}^{(n)}, \mathbf{A}^{(n)})$ rooted at a uniform point is equal to $ \mathfrak{t}$ is

\begin{align*}&\sum_{\substack{i_1, \dots, i_k \\  \text{all different}}} \mathbb{P} \left( \forall 1 \leq j \leq k,\ (I_{i_j}^{(n)},A_{i_j}^{(n)}) = (d_{v_j}^{ \mathrm{in}}, d_{v_j}^{ \mathrm{out}})\right)\underbrace{  \frac{1}{n} }_{ \text{root label}} \underbrace{ \frac{ (S_{ \mathrm{in}}-k)!}{S_{ \mathrm{in}}!}}_{\substack{\text{connection of} \\ \text{matched legs}}}  \underbrace{ \prod_{j=0}^{S_{out}- k-1} \left(1 - \frac{A}{S_{ \mathrm{in}} - k - j }\right)}_{ \text{unmatched legs}} \\
&= \left( \frac{1}{n} + o \left( \frac{1}{n}\right)\right)^k \sum_{\substack{i_1, \dots, i_k \\\text{all different}}} \mathbb{P} \left( \forall 1 \leq j \leq k,\ (I_{i_j}^{(n)},A_{i_j}^{(n)}) = (d_{v_j}^{ \mathrm{in}}, d_{v_j}^{ \mathrm{out}})\right) \cdot \prod_{j= 0}^{A-1} \left( \frac{S_{ \mathrm{in}} - S_{ \mathrm{out}} - j }{S_{ \mathrm{in}} - k - j}\right). \\
\end{align*}
Using now the fact that we are under \eqref{eq:hypo-convergence}, we obtain that the (unconditionned) probability that $ \mathrm{FCM}( \mathbf{I}^{(n)}, \mathbf{A}^{(n)})$ rooted at a uniform point is equal to $ \mathfrak{t}$  converges as $n \to \infty$ towards $$ \left(\prod_{ v \in \mathbf{t}} \nu_{d^{ \mathrm{in}}_{v}} \mu_{(d^{ \mathrm{in}}_{v} ),d^{ \mathrm{out}}_{v}} \right)\cdot (1- \mathbb{E}_{ \nu} [m])^A = \mathbb{P} \left( \mathcal{T} = \mathfrak{t} \right).$$
Since the tree $ \mathcal{T}$ is almost surely finite when $\E_{\overline{\nu}} [m] \leq 1$ and $ \Theta \geq 0$, the local limit of $ \mathrm{FCM}( \mathbf{I}^{(n)}, \mathbf{A}^{(n)})$ has the same law as $ \mathcal{T}$ in this case, which concludes the proof.
\end{proof}

We have now all the tools to prove the characterization of the phase transition given in Theorem~\ref{thm:coupling:curienhenard}. 
\proof[Proof of Theorem \ref{thm:coupling:curienhenard}] We assume $\E_{\overline{\nu}} [m] \leq 1$ and $\E_{{\nu}} [m] \leq 1$. First, thanks to Corollary \ref{cor:phasetransition}, it only remains to show that $ \Theta \geq 0$  if and only if $ \mathbb{P}\left(\exists \text{ infinite parked cluster in } T_{ \infty} ( \nu)\right) = 0$. Recall that by Proposition \ref{prop:couplingmapping} (for the mapping case) and Proposition \ref{prop:couplingtree} (for the tree case), the nearly parked components are the same as the weakly connected components of the frozen configuration model. When $ \Theta \geq 0$, the local limit of the frozen configuration model is almost surely finite by Proposition \ref{prop:souscrit}, and thus the probability that the root $ \rho_ \infty$ of $ \mathcal{T}_{ \infty} ( \nu)$ belongs to an infinite cluster of parked cars is equal to $0$. But if there exists a vertex which has a positive probability to be in an infinite cluster, then the same must be true for the root, thus
$ \mathbb{P}\left(\exists \text{ infinite parked cluster in } T_{ \infty} ( \nu)\right) = 0$. 

Conversely, suppose now that $ \mathbb{P} \left( \exists \text{ infinite parked cluster}\right) = 0$. Then the local limit of parked clusters are finite trees. It follows that the local limit of the frozen configuration model is also supported by finite trees, which are necessarily white (since the blue components have a cycle). Moreover, in white trees all out-edges are matched thus the proportion of unmatched in-legs is $ 1- \mathbb{E}_ \nu[m]$ and their are uniformly spread among the in-legs by construction. Thus, the local limit of the frozen configuration model is the Bienaym\'e--Galton--Watson tree that we described in the proof of Proposition \ref{prop:souscrit} and this tree is subcritical implying $ \Theta \geq 0$. This concludes the proof.

\endproof

\section{Conclusion and Perspectives}

The main input of this paper is the coupling between the parking process and an oriented configuration model. It allows us to study precisely the size of the connected components, in particular in the critical regime. 
More precisely, when the limit law $ \lambda$ that appears in \eqref{eq:hypo-convergence} is bounded or better when the third moment of $ \nu$ is finite and the third moments of the laws $ \mu_{(k)}$ are bounded, the behavior of the configuration model is close to that of the Erd\H{o}s--Rényi random graph. Thus, we expect that the behavior of the frozen version of both graphs is similar, and as described in \cite{contat2021parking}, we expect to have parked components of size of order $n^{2/3}$ at the critical point, and a flux of outgoing cars of order $n^{1/3}$.

However, different universality classes appear for the (unoriented) critical configuration model when the law of the degree sequence has a finite second moment but an infinite third moment and \emph{heavy tails}. For example, if the degree distribution is such that the probability that a vertex has degree $k$ is asymptotic to $c \cdot k^{- \gamma}$ for some constant $c>$ and some $ \gamma \in (3,4)$, then the sizes of connected components are of order $n^{ (\gamma -2)/( \gamma - 1)}$, see for example  \cite{conchon2023stable,dhara2020heavy,joseph2014component}. 

In our case, we consider an (frozen) oriented version of the configuration model but are interested in the weakly connected components, which is not standard in the literature. We though believe that the same phenomenon appears. In particular, when the offspring distribution of the tree or the car arrivals distributions have an infinite third moments and a tail of order $c \cdot k^{- \gamma}$ for some constant $c>$ and some $ \gamma \in (3,4)$ when $k$ goes to infinity, the size of the components should be of order  $n^{ (\gamma -2)/( \gamma - 1)}$ at criticality. 

 These exponents already appear in the enumeration of plane fully parked trees performed by Chen \cite{chen2021enumeration}. Indeed, he found that when the car arrival distribution has a tail of order $c \cdot k^{- \gamma}$ for some constant $c>$ and some $ \gamma \in (3,4)$ when $k$ goes to infinity, then the total weight of the fully parked trees of size $n$ is of order $n^{\frac{ \gamma -1}{ \gamma -2} - 1}$, which suggests also that the critical components are of order $n^{ (\gamma -2)/( \gamma - 1)}$ at criticality.

\bibliographystyle{siam}
\bibliography{/Users/contat/Dropbox/Articles/biblio}

\begin{thebibliography}{10}

\bibitem{addario2023foata}
{\sc L.~Addario-Berry, A.~Blanc-Renaudie, S.~Donderwinkel, M.~Maazoun, and
  J.~B. Martin}, {\em The {F}oata--{F}uchs proof of {C}ayley's formula, and its
  probabilistic uses}, Electronic Communications in Probability, 28 (2023),
  pp.~1--13.

\bibitem{aldous1991asymptotic}
{\sc D.~Aldous}, {\em Asymptotic fringe distributions for general families of
  random trees}, The Annals of Applied Probability,  (1991), pp.~228--266.

\bibitem{aldous2022parking}
{\sc D.~Aldous, A.~Contat, N.~Curien, and O.~H{\'e}nard}, {\em Parking on the
  infinite binary tree}, Probability Theory and Related Fields,  (2023),
  pp.~1--24.

\bibitem{chen2021enumeration}
{\sc L.~Chen}, {\em Enumeration of fully parked trees}, arXiv preprint
  arXiv:2103.15770,  (2021).

\bibitem{chen2021parking}
{\sc Q.~Chen and C.~Goldschmidt}, {\em Parking on a random rooted plane tree},
  Bernoulli, 27 (2021), pp.~93--106.

\bibitem{conchon2023stable}
{\sc G.~Conchon-Kerjan and C.~Goldschmidt}, {\em The stable graph: the metric
  space scaling limit of a critical random graph with iid power-law degrees},
  The Annals of Probability, 51 (2023), pp.~1--69.

\bibitem{contat2020sharpness}
{\sc A.~Contat}, {\em Sharpness of the phase transition for parking on random
  trees}, Random Structures \& Algorithms, 61 (2022), pp.~84--100.

\bibitem{contat2022surprising}
\leavevmode\vrule height 2pt depth -1.6pt width 23pt, {\em Surprising
  identities for the greedy independent set on {C}ayley trees}, Journal of
  Applied Probability, 59 (2022), pp.~1042--1058.

\bibitem{contat2023parking}
\leavevmode\vrule height 2pt depth -1.6pt width 23pt, {\em Parking on random
  trees}, PhD thesis, Universit{\'e} Paris-Saclay, 2023.

\bibitem{contat2021parking}
{\sc A.~Contat and N.~Curien}, {\em Parking on {C}ayley trees and frozen
  {E}rd{\H{o}}s--{R}{\'e}nyi}, The Annals of Probability, 51 (2023),
  pp.~1993--2055.

\bibitem{curien2022phase}
{\sc N.~Curien and O.~H{\'e}nard}, {\em The phase transition for parking on
  {G}alton--{W}atson trees}, Discrete Analysis,  (2022).

\bibitem{dhara2020heavy}
{\sc S.~Dhara, R.~van~der Hofstad, J.~S. van Leeuwaarden, and S.~Sen}, {\em
  Heavy-tailed configuration models at criticality}, in Annales de l'Institut
  Henri Poincar{\'e}-Probabilit{\'e}s et Statistiques, vol.~56, 2020,
  pp.~1515--1558.

\bibitem{GP19}
{\sc C.~Goldschmidt and M.~Przykucki}, {\em Parking on a random tree},
  Combinatorics, Probability and Computing, 28 (2019), pp.~23--45.

\bibitem{harary1964number}
{\sc F.~Harary, G.~Prins, and W.~Tutte}, {\em The number of plane trees},
  Indag. Math, 26 (1964), pp.~319--329.

\bibitem{Jan12b}
{\sc S.~Janson}, {\em Simply generated trees, conditioned {G}alton--{W}atson
  trees, random allocations and condensation.}, Probability Surveys, 9 (2012),
  pp.~103--252.

\bibitem{JO18}
{\sc O.~D. Jones}, {\em Runoff on rooted trees}, 2018.

\bibitem{joseph2014component}
{\sc A.~Joseph}, {\em The component sizes of a critical random graph with given
  degree sequence}, Annals of Applied Probability, 24 (2014), pp.~2560--2594.

\bibitem{konheim1966occupancy}
{\sc A.~G. Konheim and B.~Weiss}, {\em An occupancy discipline and
  applications}, SIAM Journal on Applied Mathematics, 14 (1966),
  pp.~1266--1274.

\bibitem{kryven2016emergence}
{\sc I.~Kryven}, {\em Emergence of the giant weak component in directed random
  graphs with arbitrary degree distributions}, Physical Review E, 94 (2016),
  p.~012315.

\bibitem{LaP16}
{\sc M.-L. Lackner and A.~Panholzer}, {\em Parking functions for mappings},
  Journal of Combinatorial Theory, Series A, 142 (2016), pp.~1 -- 28.

\end{thebibliography}
\end{document}